\documentclass[11pt]{amsart}
\usepackage[utf8]{inputenc}
\usepackage{amsmath}
\usepackage{mathtools}
\usepackage{thmtools}
\usepackage[colorlinks]{hyperref}
\usepackage{amsthm}
\usepackage{amssymb}
\usepackage[left=1.25in, right=1.25in, top=1.25in, bottom=1.25in]{geometry}
\usepackage{tikz}
\usepackage{microtype}
\usepackage{booktabs}
\usepackage{xurl}
\usepackage{cleveref}

\makeatletter
\def\@makefnmark{%
  \leavevmode
  \raise.9ex\hbox{\check@mathfonts\fontsize\sf@size\z@\normalfont\@thefnmark}%
}
\makeatother

\numberwithin{equation}{section}

\newtheorem{theorem}{Theorem}

\newtheorem{conjecture}[theorem]{Conjecture}

\theoremstyle{definition}
\newtheorem{example}[theorem]{Example}

\theoremstyle{remark}

\crefname{conjecture}{Conjecture}{Conjectures}
\Crefname{conjecture}{Conjecture}{Conjectures}

\title{A counterexample to a log-concavity conjecture of Brenti}
\author{Christian Gaetz}
\address[Gaetz]{Department of Mathematics, University of California, Berkeley, CA.}
\email{{\href{mailto:gaetz@berkeley.edu}{{gaetz@berkeley.edu}}}}
\date{\today}
\begin{document}

\begin{abstract}
This note records a counterexample to Brenti's conjecture (\emph{Discrete Math.}, 1998) that the nonzero coefficients of $\widetilde{R}$-polynomials of the symmetric group form a log-concave sequence.
\end{abstract}


\maketitle

\section{Introduction}

The log-concavity and unimodality of various combinatorial sequences have in recent years been of great interest, in large part due to a series of spectacular positive results for several invariants of matroids \cite{Adiprasito-Huh-Katz, Anari-Liu-OveisGharan-Vinzant-III, Branden-Huh}. 

On the other hand, Larson \cite{Larson-counterexamples} has used AI tools to disprove other conjectures in this area. Cheng--Liu \cite{Cheng-Liu} have likewise used AI tools to disprove unimodality conjectures for \emph{Kazhdan--Lusztig polynomials} \cite{Braden-Huh-Matherne-Proudfoot-Wang,Elias-Proudfoot-Wakefield} of matroids.

In this note we use AI tools to disprove a log-concavity conjecture of Brenti \cite{Brenti-KL-combinatorial} in the (more classical) Kazhdan--Lusztig theory of Coxeter groups.

\section{On log-concavity for \texorpdfstring{$\widetilde{R}$}{R-tilde}-polynomials}
Let $\leq$ denote Bruhat order on a Coxeter group $W$. For $u,v \in W$, Kazhdan and Lusztig defined the $R$-polynomial $R_{u,v}(q) \in \mathbb{Z}[q]$, which is nonzero if and only if $u \leq v$ \cite{Kazhdan-Lusztig}. The $\widetilde{R}$-polynomials $\{\widetilde{R}_{u,v}(q)\}_{u,v \in W} \subset \mathbb{N}[q]$ are defined uniquely by the condition that 
\[
q^{\ell(u,v)/2} \widetilde{R}_{u,v}(q^{1/2}-q^{-1/2})=R_{u,v}(q),
\]
where $\ell(u,v) \coloneqq \ell(v)-\ell(u)$ is the relative length function on $W$. Combinatorial and geometric interpretations of their coefficients were given by Deodhar \cite{Deodhar-Bruhat-orderings-I} and by Dyer \cite{Dyer-thesis,Dyer-shellings}. 

All powers of $q$ appearing in a given $\widetilde{R}$-polynomial have the same parity as $\ell(u,v)$, so there exist polynomials $\{Q_{u,v}(q)\}_{u,v \in W} \subset \mathbb{N}[q]$ such that 
\[
\widetilde{R}_{u,v}(q) = \begin{cases} Q_{u,v}(q^2) & \text{if $\ell(u,v)$ is even,} \\ qQ_{u,v}(q^2) & \text{if $\ell(u,v)$ is odd.} \end{cases}
\]
Moreover, the sequence $(a_0,a_1,\ldots,a_{\lfloor \ell(u,v)/2 \rfloor})$ of coefficients of $Q_{u,v}(q)=\sum_i a_iq^i$ has no internal zeros. Recall that this sequence, and the polynomial $Q_{u,v}$ itself, is called \emph{log-concave} if $a_i^2 \geq a_{i-1}a_{i+1}$ for all $1 \leq i \leq \lfloor \ell(u,v)/2 \rfloor-1$. Brenti \cite{Brenti-KL-combinatorial} conjectured that this property is always enjoyed by the $Q_{u,v}$ polynomials for the symmetric group.

\begin{conjecture}[Brenti {\cite[Conj.~7.1]{Brenti-KL-combinatorial}}]
\label{conj:typeA}
Let $W$ be the symmetric group $S_n$. Then $Q_{u,v}$ is log-concave for all permutations $u \leq v$.
\end{conjecture}

More recently, Brenti extended the conjecture to all finite Coxeter groups $W$ \cite{Brenti-open-problems}.

\begin{conjecture}[Brenti {\cite[Conj.~2.4]{Brenti-open-problems}}]
\label{conj:finite-W}
Let $W$ be a finite Coxeter group. Then $Q_{u,v}$ is log-concave for all $u \leq v$ in $W$.
\end{conjecture}

\begin{theorem}
\label{thm:counterexample}
\Cref{conj:typeA,conj:finite-W} are false.
\end{theorem}

\Cref{thm:counterexample} is demonstrated by means of explicit counterexamples. \Cref{conj:finite-W} is easier to falsify, as there exist several\footnote{While writing this article, the author learned of another recently discovered counterexample in type $H_4$ \cite{Paul-C-LinkedIn}.} counterexamples\footnote{After posting the first version of this article, the author learned of a third previously discovered $H_4$ counterexample \cite{Eberhardt}.} in type $H_4$.

\begin{example}
\label{ex:H4}
Let $W$ be of type $H_4$ with $m(s_1,s_2)=5$, $m(s_2,s_3)=m(s_3,s_4)=3$, and other pairs of simple generators commuting. Let $u < v$ be given by the reduced words:
\begin{align*}
    u &= s_2 s_3 s_2 s_1 s_2 s_1 s_4 s_3 s_2 s_1 s_2 s_1 s_3 s_2 s_1 s_2 s_3 s_4 s_3 s_2 s_1, \\
    v &= s_1 s_2 s_1 s_2 s_3 s_2 s_1 s_2 s_1 s_3 s_2 s_4 s_3 s_2 s_1 s_2 s_1 s_3 s_2 s_1 s_2 \\
      &\qquad s_3 s_4 s_3 s_2 s_1 s_2 s_1 s_3 s_2 s_1 s_2 s_4 s_3 s_2 s_1 s_2 s_3 s_4.
\end{align*}
Then
\[
Q_{u,v}(q) = q^9 + 14q^8 +78q^7+220q^6+326q^5+234q^4+67q^3+8q^2+q,
\]
and we see that $a_2^2=8^2=64 < 67 = a_1a_3$, contradicting \Cref{conj:finite-W}.
\end{example}

The author did not find a counterexample to \Cref{conj:typeA} in $S_8$, the largest symmetric group he was able to search exhaustively. However, the following example in $S_{14}$ was constructed.

\begin{example}
\label{ex:S14}
Let $W$ be the symmetric group $S_{14}$ and let $u<v$ be the permutations whose one-line notations are:
\begin{align*}
    u &= [1,2,5,7,9,3,11,4,6,12,13,8,10,14], \\
    v &= [8,4,6,12,10,13,1,14,2,11,3,7,5,9].
\end{align*}
Then 
\[
Q_{u,v}(q) = q^{14} + 16q^{13} + 101q^{12} + 333q^{11} + 630q^{10} + 695q^{9} + 425q^{8} + 123q^{7} + 11q^6 + q^5,
\]
and we see that $a_{6}^2=11^2=121<123=a_5a_7$, contradicting \Cref{conj:typeA}.
\end{example}

\section*{Statement on AI use}
\Cref{ex:H4} was found by Fable 5. ChatGPT 5.6 Pro then found a counterexample in $S_{14}$ after being instructed to try to replicate features of the $H_4$ example. The original $S_{14}$ counterexample was simplified by applying the descent recurrence for $\widetilde{R}$-polynomials (shortening $u$ and $v$).

Both counterexamples were verified by the author using the \texttt{coxeter3} package in SageMath \cite{sagemath}.

The paper was written by the author, with minor editing carried out by AI tools.

\section*{Acknowledgments}
I am grateful to Jens Eberhardt for alerting me to the type $H_4$ counterexample he discovered. The author was partially supported by NSF grant DMS-2452032 and by a travel grant from the Simons Foundation. 

\bibliographystyle{halpha-abbrv}
\bibliography{references}

\end{document}